\documentclass[amstex,12pt,russian,amssymb]{article}

\usepackage{mathtext}
\usepackage[cp1251]{inputenc}
\usepackage[T2A]{fontenc}
\usepackage[russian]{babel}
\usepackage[dvips]{graphicx}
\usepackage{amsmath}
\usepackage{amssymb}
\usepackage{amsxtra}
\usepackage{latexsym}
\usepackage{ifthen}

\renewcommand{\baselinestretch}{1.3}

\begin{document}

\def\Xint#1{\mathchoice
   {\XXint\displaystyle\textstyle{#1}}%
   {\XXint\textstyle\scriptstyle{#1}}%
   {\XXint\scriptstyle\scriptscriptstyle{#1}}%
   {\XXint\scriptscriptstyle\scriptscriptstyle{#1}}%
   \!\int}
\def\XXint#1#2#3{{\setbox0=\hbox{$#1{#2#3}{\int}$}
     \vcenter{\hbox{$#2#3$}}\kern-.5\wd0}}
\def\dashint{\Xint-}

\renewcommand{\abstractname}{}

\noindent
\renewcommand{\baselinestretch}{1.1}
\normalsize \large \noindent

\title{Асимптотическое поведение на бесконечности //
  решений уравнения Бельтрами}

\author{Е.С. Афанасьева, Р.Р. Салимов}

\medskip

 {УДК 517.5}

{\bf Е.С. Афанасьева}  {(Ин-т прикладной математики и механики НАН
Украины, Славянск)}

{\bf  Р.Р. Салимов}  {(Ин-т  математики  НАН
Украины, Киев)}
\medskip

{\bf  Асимптотическое поведение на бесконечности решений  уравнения Бельтрами.}

\begin{abstract} В данной статье исследуется асимптотическое поведение на бесконечности гомеоморфных
решений уравнения Бельтрами при различных условиях на дилатацию.

\end{abstract}

{\bf 1. Введение.}  Пусть $D$ --  область в комплексной плоскости
${\Bbb C}$, т.е. связное открытое подмножество ${\Bbb C}$ и пусть
$\mu(z): D\rightarrow \mathbb{C}$ -- измеримая функция с
$|\mu(z)|<1$ п.в. (почти всюду) в $D$. {\it Уравнением Бельтрами}
называется уравнение вида \begin{equation}\label{eqBeltrami} f_{\bar
z}=\mu(z)\,f_z,\end{equation} где $f_{\bar
z}=\overline{\partial}f=(f_x+if_y)/2$, $f_{z}=\partial
f=(f_x-if_y)/2$, $z=x+iy$,
 $f_x$ и $f_y$ частные производные отображения $f$ по $x$ и $y$, соответственно.
 Функция $\mu$ называется {\it комплексным коэффициентом}, а
\begin{equation}\label{eqKPRS1.1}K_{\mu}(z)=\frac{1+|\mu(z)|}{1-|\mu(z)|}\end{equation} --
{\it  дилатационным отношением} уравнения (\ref{eqBeltrami}).
Уравнение Бельтрами (\ref{eqBeltrami})
называется {\it
вырожденным}, если ${\rm ess}\,{\rm sup}\,K_{\mu}(z)=\infty$.

Существование гомеоморфного $W^{1,1}_{\rm loc}$ решения было недавно
установлено для многих вырожденных уравнениях Бельтрами,
 см., напр., соответствующие ссылки в монографиях
\cite{GRSY^*}, \cite{MRSY}.

\medskip

{\bf 2. Предварительные сведения.} Напомним некоторые определения. Борелева функция
$\rho:\mathbb{C}\to[0,\infty]$ называется {\it допустимой} для
семейства кривых $\Gamma$ в $\mathbb{C}$, пишут $\rho\in{\rm
adm}\,\Gamma$, если
\begin{equation}\int\limits_{\gamma}\rho(z)\,|dz|\geqslant 1 \end{equation}
для всех $\gamma\in\Gamma$. Тогда
{\it модулем} семейства кривых $\Gamma$ называется  величина

\begin{equation}
\mathcal{M}(\Gamma)=\inf\limits_{\rho\in{\rm adm}\,\Gamma}
\int\limits_{\mathbb{C}}\rho^{2}(z)\,dm(z)
\end{equation}
Здесь $m$ обозначает меру Лебега в $\mathbb{C}$.

Следуя работе \cite{MRV}, пару $\mathcal{E}=(A,C)$, где $A\subset\mathbb{C}$
-- открытое множество и $C$ -- непустое компактное множество,
содержащееся в $A$, называем {\it конденсатором}. Конденсатор $\mathcal{E}$
называется  {\it кольцевым конденсатором}, если $G=A\setminus C$ --
кольцо, т.е., если $G$ -- область, дополнение которой
$\overline{\mathbb{C}}\setminus G$ состоит в точности из двух
компонент.  Говорят
также, что конденсатор $\mathcal{E}=(A,C)$ лежит в области $D$, если $A\subset
D$. Очевидно, что если $f:D\to\mathbb{C}$ -- непрерывное, открытое
отображение и $\mathcal{E}=(A,C)$ -- конденсатор в $D$, то $(fA,fC)$ также
конденсатор в $fD$. Далее $f\mathcal{E}=(fA,fC)$.

Функция $u:A\to \mathbb{R}$ {\it абсолютно непрерывна на прямой}, имеющей непустое пересечение с $A$, если она абсолютно непрерывна на любом отрезке этой прямой, заключенном в $A$. Функция $u:A\to \mathbb{R}$ принадлежит классу ${\rm ACL}$ ({\it абсолютно непрерывна на почти всех прямых}), если она абсолютно непрерывна на почти всех прямых, параллельных любой координатной оси.

Обозначим через $C_0(A)$  множество
непрерывных функций $u:A\to\mathbb{R}$ с компактным носителем,
$W_0(\mathcal{E})=W_0(A,C)$ -- семейство неотрицательных функций
$u:A\to\mathbb{R}$ таких, что 1) $u\in C_0(A)$, 2)
$u(z)\geqslant 1$ для $z\in C$ и 3) $u$ принадлежит классу ${\rm
ACL}$. Также обозначим

\begin{equation}
\vert\nabla
u\vert=\sqrt{{\,{\left(\frac{\partial u}{\partial x}\right)}^2+\left(\frac{\partial u}{\partial y}\right)}^2
}.
\end{equation}

Величину
\begin{equation}
{\rm cap}\,\mathcal{E}={\rm cap}\,(A,C)=\inf\limits_{u\in W_0(\mathcal{E})}\,
\int\limits_{A}\,\vert\nabla u\vert^2\,dm(z)
\end{equation}
называют
{\it ёмкостью} конденсатора $\mathcal{E}$.

Пусть $D$ -- область в  ${\Bbb C}$.  $E\,,F\,\subseteq\, D\,-$
произвольные множества. Обозначим через $\Delta(E,F;D)$ семейство
всех кривых $\gamma:[a,b]\,\rightarrow\,{\Bbb C}\,,$ которые
соединяют $E$ и $F$ в $D\,,$ т.е. $\gamma(a)\,\in\,E\,,\gamma(b)
\,\in\,F$ и $\gamma(t)\,\in\,D\,$ при $a\,<\,t\,<\,b\,.$

В дальнейшем  мы
будем использовать равенство
\begin{equation}\label{EMC}
{\rm cap}\,\mathcal{E}=\mathcal{M}(\Delta(\partial A,\partial C; A\setminus C)),\ \
 \end{equation} см. теорему  1 в \cite{Sh}.

Известно, что \begin{equation}\label{maz} {\rm
cap}\,\mathcal{E}\geqslant \frac{4\pi}{\ln \frac{m(A)}{m(C)}}
\end{equation}
см., напр., неравенство (8.9)  в \cite{MazIS}.

Напомним следующие термины.  Пусть $d_0\,=\,\rm {dist}\,(z_0\,,\partial D)$ и
пусть $Q:D\rightarrow\,[0\,,\infty]\,$ -- измеримая по Лебегу функция.
Положим

\begin{equation}\label{2.10}
\mathbb{A}(z_0,r_1,r_2) = \{ z\,\in\,{\Bbb C} : r_1<|z-z_0|<r_2\}\ ,
\end{equation}

\begin{equation}\label{2.11}
S_i\,=\,S(z_0,r_i) = \{ z\,\in\,{\Bbb C} : |z-z_0|=r_i\}\,\,,\ \ \
i=1,2.
\end{equation}

 Будем говорить, что гомеоморфизм $f:D\to {\Bbb C}$ является
{\it кольцевым $Q$-го\-ме\-о\-мор\-физ\-мом
в точке $z_0\,\in\,D,$} если соотношение

\begin{equation}\label{defring}
\mathcal{M}\,\left(\Delta\left(fS_1,fS_2, fD\right)\right)\,\ \leqslant
\int\limits_{\mathbb{A}} Q(z)\cdot \eta^{2}(|z-z_0|)\ dm(z)
\end{equation}
выполнено для любого кольца $\mathbb{A}=\mathbb{A}( z_0,r_1,r_2),$\,\, $0<r_1<r_2<
d_0$ и для каждой измеримой функции $\eta : (r_1,r_2)\to [0,\infty
]\,,$ такой, что
\begin{equation}\label{norma}
\int\limits_{r_1}^{r_2}\eta(r)\ dr\ \geqslant\ 1\,.
\end{equation}
Говорят, что гомеоморфизм $f:D\to {\Bbb C}$ является {\it
кольцевым $Q$-го\-ме\-о\-мор\-физ\-мом  в
области $D$}, если условие $(\ref{defring})$ выполнено для всех точек
$z_0\,\in\,D\,.$

Следующее  утверждение можно найти в работе \cite{KSS}, теорема 5.1.

\medskip

\textbf{Теорема 1.} \rm {\it \rm {\it
 Пусть $D$ и  $D'$ -- области в $\mathbb{C}$, и   $f:D\rightarrow D'$ --
гомеоморфное   решение уравнения Бельтрами (\ref{eqBeltrami}) класса $W^{1,1}_{\rm
loc}$ и $K_{\mu}\in L^1_{\rm loc}(D)$. Тогда  $f$ является кольцевым  $Q$-гомеоморфизмом в каждой точке $z_0\in D$    с
$Q(z)=K_{\mu}(z)$.}

\medskip

{\bf 3. Поведение на бесконечности.}
Ассимтотическое поведение на бесконечности  кольцевых $Q$-гомеоморфизмов при оптимальных условиях исследовалось в работе \cite{SS2}.  Пусть $r_0$ -- произвольное фиксированное положительное число. Для гомеоморфизма $f:\mathbb{C}\to \mathbb{C}$    полагаем

\begin{equation}
M(R,f)=\max\limits_{|z-z_0|=R}\, |f(z)-f(z_0)|\,.
\end{equation}

\medskip

\textbf{Лемма 2.} \rm {\it \rm {\it
Пусть $f:\mathbb{C}\to \mathbb{C}$ -- гомеоморфное     решение уравнения Бельтрами (\ref{eqBeltrami}) класса $W^{1,1}_{\rm
loc}$. Если  $K_{\mu}\in L^1_{\rm loc}(\mathbb{C})$, тогда
\begin{equation}\label{SI}
\mathcal{M}\,\left(\Delta\left(fS_1,fS_2, f\mathbb{A}\right)\right)\,\ \leqslant \Lambda(R)\,,
\end{equation}
где $S_1=S(z_0, r_0), \, S_1=S(z_0, R) $, $\mathbb{A}=\mathbb{A}(z_0, r_0, R)$
и
\begin{equation}\label{90}
\Lambda(R)= \left(\int\limits_{r_0}^R  \psi(t) \, dt\right)^{-2}\,\cdot\int\limits_{\mathbb{A}}K_{\mu}(z)\,\psi^2(|z-z_0|) \, dm(z)
\end{equation}
для любой измеримой  (по Лебегу) функции  $\psi:[0,\infty]\to [0,\infty]$ такой, что
\begin{equation}\label{I}
0<\int\limits_{r_0}^R  \psi(t) \, dt<\infty \ \ \ \ \ \ \ \ \forall \ R>r_0\,.
 \end{equation}}

\medskip

{\it Доказательство.} Пусть
$\mathbb{A}=\mathbb{A}(z_0,r_0, R)$ с $0<r_0<R$.
Рассмотрим измеримую функцию
\begin{equation}
\eta(t)\,=\,\left
\{\begin{array}{rr} \frac{\psi(t)}{\int\limits_{r_0}^{R}\, \psi(t)\, dt}, &  \ t\in (r_0,R) \\
0, & \ t\in \Bbb{R}\setminus (r_0,R).
\end{array}\right.
\end{equation}
Отметим, что функция $\eta(t)$ удовлетворяет условию (\ref{norma}). Тогда из  теоремы 1 вытекает оценка (\ref{SI}).

\medskip

\textbf{Лемма 3.} \rm {\it \rm {\it
Пусть $\mu :\mathbb{C}\to \mathbb{C}$ -- измеримая функция с $|\mu(z)|<1$ п.в. такая, что   $K_{\mu}\in L^1_{\rm loc}(\mathbb{C})$.
Тогда уравнение
(\ref{eqBeltrami}) не имеет гомеоморфного   решения $f:\mathbb{C}\to \mathbb{C}$ класса Соболева $W_{\rm loc}^{1,1}$  с асимптотикой \begin{equation}\label{as1}
\liminf\limits_{R\to \infty} M(z_0,f,R) \, e^{-\frac{2\pi}{\Lambda(R)}}=0\,.
\end{equation}
}

\medskip

{\it Доказательство.}  Предположим противное, а именно, что существует гомеоморфное решение $f:\mathbb{C}\to \mathbb{C}$ класса $W^{1,1}_{\rm
loc}$.

Рассмотрим  кольцо
$\mathbb{A}=\mathbb{A}(z_0,r_0, R)$ с $0<r_0<R$. Тогда
$\left(f B_R,f\overline{B}_{r_0}\right)$
-- кольцевой конденсатор в   $\mathbb{C}$ и, согласно (\ref{EMC}), имеем
равенство
$$ {\rm cap}\left(f B_R,f\overline{B}_{r_0}\right) =\mathcal{M}(\triangle(\partial
fB_R,\partial f B_{r_0};f\mathbb{A}))$$ а ввиду
гомеоморфности  $f,$ равенство
$$\triangle\left(\partial
fB_R,\partial
fB_{r_0};f\mathbb{A}\right)=f\left(\triangle\left(\partial
B_R,\partial
B_{r_0};\mathbb{A}\right)\right).$$

В  силу леммы 2  имеем

\begin{equation}\label{leq101} {\rm cap}\ \left(fB_R,\overline{fB}_{r_0}\right) \leqslant
\Lambda(R)\ .\end{equation}
С другой стороны,  в силу  неравенства (\ref{maz}) вытекает оценка

\begin{equation}\label{eq102} {\rm cap}\ \left(fB_R,\overline{fB}_{r_0}\right)
\geqslant  \frac{4\pi}{\ln \, \frac{m\left(fB_R\right)}{m\left(\overline{fB}_{r_0}\right)}}\,.
\end{equation}

Комбинируя   (\ref{leq101}) и (\ref{eq102}), получаем, что
\begin{equation}\label{U}
m\left(\overline{fB}_{r_0}\right)\leqslant m\left(fB_R\right) e^{-\frac{4\pi}{\Lambda(R)}}\,.
\end{equation}

Заметим, что $m\left(fB_R\right) \leqslant \pi   M^2(z_0,f,R)$, поэтому из неравенства (\ref{U}) вытекает следующая оценка

\begin{equation}\label{U1}
\sqrt{\frac{m\left(\overline{fB}_{r_0}\,\right)}{\pi}}\leqslant M(z_0,f,R) \, e^{-\frac{2\pi}{\Lambda(R)}}\,.
\end{equation}

Очевидно,  $M_0=\sqrt{\frac{m\left(\overline{fB}_{r_0}\,\right)}{\pi}}>0$  и не зависит от $R$.
Переходя к  нижнему  пределу при
$R\to \infty$ и учитывая условие (\ref{as1}),  получаем $m(fB_{r_0})=0$, что противоречит гомеоморфности отображения $f$.

\medskip

\textbf{Лемма 4.} \rm {\it \rm {\it  Пусть $\mu :\mathbb{C}\to \mathbb{C}$ -- измеримая функция с $|\mu(z)|<1$ п.в.  и   $K_{\mu}\in L^1_{\rm loc}(\mathbb{C})$. Предположим, что существует  неотрицательная измеримая   (по Лебегу) функция такая, что
\begin{equation}\label{I}
0<I(R)=\int\limits_{r_0}^R  \psi(t) \, dt<\infty \ \ \ \ \ \ \ \ \forall\  R\ >r_0\,,\end{equation}
и при $R\to \infty$
\begin{equation}\label{SI1*}
\int\limits_{\mathbb{A}(z_0,r_0, R)} K_{\mu}(z) \psi^2(|z-z_0|)\,dm(z)\leqslant c\cdot I^p(R)\,,
\end{equation}
где $p\leqslant 2$. Тогда уравнение
(\ref{eqBeltrami}) не имеет гомеоморфного   решения $f:\mathbb{C}\to \mathbb{C}$ класса Соболева $W_{\rm loc}^{1,1}$ с асимптотикой
\begin{equation}\label{SI1**}
\liminf\limits_{R\to \infty} M(z_0,f,R)\exp\left(- \frac{2\pi}{c}\,I^{2-p}(R)\right)=0\,.
\end{equation}}

\medskip

{\it Доказательство.} Ввиду   условия (\ref{SI1*}) получаем

 \begin{equation}
\Lambda(R)=\frac{\int\limits_{\mathbb{A}(z_0,r_0,R)}K_{\mu}(z)\,\psi^2(|z-z_0|) \, dm(z)}{\left(\int\limits_{r_0}^R  \psi(t) \, dt\right)^2}\leqslant c I^{p-2}(R)\,.
\end{equation}
Отсюда следует, что
\begin{equation}
\exp\left[-\frac{2\pi}{\Lambda(R)}\right]\leqslant \exp{\left[-{\frac{2\pi}{c}}\,I^{2-p}(R)\right]^{}}\,.
\end{equation}
Таким образом, заключение леммы 4 вытекает из леммы 3.

\bigskip

В дальнейшем, для целых $k\geqslant 0$ полагаем
\begin{equation} e_0=1,\  \ \ \  e_1=e, \  \ \ \  e_2=e^e, \  \ldots ,\  e_{k+1}=\exp e_{k}\,   \end{equation}
и
\begin{equation}  \ln_0 t=t,\  \ \ \  \ln_1 t=\ln t, \  \ \ \  \ln_2 t= \ln\ln t, \  \ldots ,\  \ln_{k+1} t= \ln \ln_{k} t \,. \end{equation}

\medskip

\textbf{Лемма 5.} \rm {\it \rm {\it При $R>e_{N}$  справедливо равенство
\begin{equation}
\int\limits_{e_{N}}^R\frac{dt }{\prod\limits_{k=0}^N \ln_k t}= \ln_{N+1} R\,.
\end{equation}}

\medskip

{\it Доказательство.}  Действительно, выполнив замену переменной  $s=\ln_N t$,  получим заявленное утверждение:
\begin{equation}
\int\limits_{e_{N}}^R\frac{dt }{\prod\limits_{k=0}^N \ln_k t}= \int\limits_{1}^{\ln_N R}  \frac{ds}{s}=\ln \ln_N R= \ln_{N+1} R\,.
\end{equation}

Выбирая в лемме 4 $\psi(t)=\frac{1}{\prod\limits_{k=0}^N \ln_k t},\ \ r_0=e_{N}$ и  $p=1$, приходим к следующему утверждению.

\medskip

\textbf{Теорема 2.} \rm {\it \rm {\it  Пусть $\mu :\mathbb{C}\to \mathbb{C}$ -- измеримая функция с $|\mu(z)|<1$ п.в. такая, что   $K_{\mu}\in L^1_{\rm loc}(\mathbb{C})$
и
\begin{equation}\label{SI1}
\int\limits_{\mathbb{A}(z_0, e_{N}, R)} \frac{K_{\mu}(z)\,dm(z)}{\left(\prod\limits_{k=0}^N \ln_k|z-z_0|\right)^2 }\leqslant c\cdot \ln_{N+1} (R)\,  \ \ \ \ \ \ \ \ \forall\  R>e_{N} \,.
\end{equation}
Тогда уравнение
(\ref{eqBeltrami}) не имеет гомеоморфного   решения $f:\mathbb{C}\to \mathbb{C}$ класса Соболева $W_{\rm loc}^{1,1}$ с асимптотикой
 \begin{equation}\label{SI12**}
 \liminf\limits_{R\to \infty}  \frac{M(z_0,f,R)}{\ln^{\gamma}_N (R)}=0\,,
\end{equation}
где $\gamma=\frac{2\pi}{C}$.}

\medskip

Выбирая в теореме  2 $N=0$, приходим к следующему следствию.

\medskip

\textbf{Следствие 1.} \rm {\it \rm {\it
Пусть $\mu :\mathbb{C}\to \mathbb{C}$ -- измеримая функция с $|\mu(z)|<1$ п.в. такая, что   $K_{\mu}\in L^1_{\rm loc}(\mathbb{C})$
и
\begin{equation}\label{SI1}
\int\limits_{\mathbb{A}(z_0, 1, R)} \frac{K_{\mu}(z)\,dm(x)}{|z-z_0|^2}\leqslant c\cdot \ln R\,  \ \ \ \ \ \ \ \ \forall \ R>1 \,,
\end{equation}
Тогда уравнение
(\ref{eqBeltrami}) не имеет гомеоморфного   решения $f:\mathbb{C}\to \mathbb{C}$ класса Соболева $W_{\rm loc}^{1,1}$ с асимптотикой
\begin{equation}\label{SI1**}
\liminf\limits_{R\to \infty}  \frac{M(z_0,f,R)}{R^{\frac{2\pi}{c}}}=0\,.
\end{equation}
}

\medskip

\textbf{Следствие 2.} \rm {\it \rm {\it
Пусть $\mu :\mathbb{C}\to \mathbb{C}$ -- измеримая функция с $|\mu(z)|<1$ п.в. такая, что   $K_{\mu}\in L^1_{\rm loc}(\mathbb{C})$
и
\begin{equation}\label{SI1}
\frac{1}{2\pi R}\int\limits_{S(z_0, R)} K_{\mu}(z)\,|dz| \leqslant K\,  \ \ \ \ \ \ \ \ \forall \ R>1 \,.
\end{equation}
Тогда уравнение
(\ref{eqBeltrami}) не имеет гомеоморфного   решения $f:\mathbb{C}\to \mathbb{C}$ класса Соболева $W_{\rm loc}^{1,1}$ с асимптотикой
\begin{equation}\label{MRE}
\liminf\limits_{R\to \infty}  \frac{M(z_0,f,R)}{R^{1/K}}=0\,.
\end{equation}}
}

\medskip

Елена Сергеевна Афанасьева

Институт прикладной математики и механики НАН Украины, Славянск

Салимов Руслан Радикович

Институт  математики  НАН Украины, Киев

Email: salimov07@rambler.ru, ruslan623@yandex.ru,


\begin{thebibliography}{30}









\bibitem{GRSY^*}
\emph{Gutlyanskii V.,  Ryazanov V.,   Srebro U.,  Yakubov E.} The
Beltrami Equation: A Geometric Approach. Springer Advances in
Mathematics. ISBN 978-1-4614-3190-9, Due: May 31, 2012.





\bibitem{MRSY}
\emph{Martio O., Ryazanov V., Srebro U.,  Yakubov E.}  Moduli in
Modern Mapping Theory. Springer Monographs in Mathematics. -- New
York: Springer, 2009. -- 367 pp.








\bibitem{MRV}
\emph{Martio O., Rickman S., Vaisala J.} Definitions for
quasiregular mappings // Ann. Acad. Sci. Fenn. Ser. A1. -- 1969. --
Vol.~448. -- P.~1--40.

\bibitem{MazIS} V.~Maz'ya, {\em Lectures on isoperimetric and isocapacitary inequalities in the theory of Sobolev spaces}.
Contemp. Math., \textbf{338} (2003), 307–-340.




\bibitem{Sh}
\emph{Шлык В.А.} О равенстве p-емкости и p-модуля
// Сиб. мат. ж. -- 1993. -- T.~34, №~6. -- C.~216–-221.


\bibitem{KSS} \emph{ Д. А. Ковтонюк, Р. Р. Салимов, Е. А. Севостьянов }  К теории отображений классов Соболева и Орлича-Соболева  - Киев : Наук. думка, 2013. - 303 с.



\bibitem{SS2} \emph{Р.Р. Салимов., Е.С. Смоловая} О порядке роста кольцевых $Q$-гомеоморфизмов на бесконечности // Укр. мат. журн. – 2010. -- Т.62, № 6. -- С. 829--836.





\end{thebibliography}
\end{document}